\documentclass[11pt,a4paper]{article}

\usepackage[T1]{fontenc}
\usepackage{lmodern}

\usepackage{amsmath}
\usepackage{amssymb}
\usepackage{graphicx}
\usepackage{float}
\usepackage{geometry}
\usepackage{hyperref}
\usepackage{authblk}

\numberwithin{equation}{section}

\title{
Compact implicit high-resolution numerical scheme for 
multicomponent transport problems with Langmuir sorption
}

\author[1]{Dagmar \v{Z}\'akov\'a\thanks{
Corresponding author: \href{mailto:dagmar.zakova@stuba.sk}
{dagmar.zakova@stuba.sk}}}

\author[1]{Peter Frolkovi\v{c}}

\affil[1]{
Department of Mathematics and Descriptive Geometry,
Faculty of Civil Engineering,
Slovak University of Technology in Bratislava,
Bratislava, Slovakia
}

\date{}

\begin{document}

\maketitle

\begin{abstract}
We present a compact implicit high-resolution finite volume scheme for one-dimensional multicomponent transport problems with nonlinear Langmuir sorption. The method combines a second order compact implicit discretization with WENO reconstruction and a time limiter for the treatment of steep gradients and discontinuities. A special attention is given to the Courant number used in the limiter, where local characteristic Courant numbers based on the eigenvalues of the inverse retardation matrix are considered. The compact implicit schemes are compared with analogous explicit first and second order methods. Numerical experiments show second order convergence for smooth solutions and demonstrate that the implicit approach can be used with significantly larger Courant numbers. The method is also applied to a three-component displacement chromatography problem with sharp concentration fronts.
\end{abstract}

\vspace{0.5em}


\vspace{0.5em}

\vspace{1em}

\section{Introduction}
\label{sec-intorduction}

For many years, there has been a growing interest in mathematical models describing the transport of substances in porous media, where the migration of solute components is significantly influenced by their interaction with the solid matrix through sorption or adsorption processes. This class of problems plays a key role in at least two important applications: the simulation of contaminant transport in groundwater systems \cite{brusseau1995effect, kacurSolutionContaminantTransport2005, frolkovic_semi-analytical_2006, patil2014contaminant, frolkovicNumericalSimulationContaminant2016} and the modeling of preparative liquid chromatography \cite{javeed2011efficient, burger_linearly_2018, meyer2020chromatech, rasheed2024simulation, zakovafrolkovic2025tatra}. 

In this paper, we focus our attention purely on the underlying mathematical model, clearly motivated by these two important engineering applications. Moreover, we concentrate only on the advective part of the transport and skip the treatment of diffusion-dispersion processes, referring to standard approaches to include them \cite{donat2024weno, frolkovicNumericalSimulationContaminant2016}. Consequently, in our study, we are dealing formally with hyperbolic partial differential equations  \cite{leveque_finite_2004,frolkovic_semi-analytical_2006} given by hyperbolic systems of first order \cite{donat2018implicit}.

For such transient problems governed by partial differential equations, to be solved numerically, two fundamental classes of time discretization can be considered. The first and often preferred approach is represented by explicit methods, which deliver numerical solutions directly without the need to solve algebraic systems of equations. However, explicit schemes are constrained by stability restrictions on the choice of the time step, such as the Courant–Friedrichs–Lewy (CFL) condition \cite{leveque_finite_2004}. 
Conversely, the second approach consists of implicit methods. Although, in general, they lead to the necessity of solving systems of algebraic equations at each time level, they offer unconditional stability, allowing the time step to be selected based solely on accuracy requirements rather than stability constraints \cite{leveque_finite_2004, zakovafrolkovivc2026esaim}.

The situation becomes significantly different when the mathematical model contains a system of differential equations with nonlinearity under the time derivative, as occurs in the presence of nonlinear accumulation terms due to, e.g., sorption isotherms. In this case, the problem to be solved is described by a system of strongly coupled transport-reaction equations, where even explicit methods are forced to solve local, nonlinear algebraic equations at each grid cell to recover physical concentrations from the conserved variables under the time derivative. Furthermore, any standard explicit discretization of each individual differential equation in a component-wise manner fails to respect the strong coupling with the remaining equations in the system.

For this class of problems, implicit methods of time discretization appear to be highly advantageous, provided that they can bypass the computationally expensive task of solving large, globally coupled, nonlinear algebraic systems. To address this challenge, we present a compact implicit numerical method \cite{frolkovic2023high, zakovafrolkovic2025tatra, zakovafrolkovivc2026esaim}, which, similar to explicit approaches, requires only the resolution of small, local systems of algebraic equations in each grid cell while fully respecting the tight coupling of the transport-reaction equations in the system.

In this paper, we develop a high-resolution, compact implicit numerical scheme that is tailored to solve transport equations coupled through fast sorption processes using popular Langmuir isotherms \cite{frolkovic_semi-analytical_2006,donat2018implicit, donat2024weno, frolkovicNumericalSimulationContaminant2016}. For smooth solutions, the method is second order accurate, and it can compete very well with analogous explicit numerical methods, as it is not restricted by the CFL condition for the choice of time steps. This is a very important advantage in cases where some of the solutes (which may not be of large interest for a particular application) are transported quickly with respect to others. If, for instance, such solutes leave the computational domain, the implicit method can enlarge the time step, which accurately approximates the slow components. This is not possible for explicit schemes, where the choice of time step must still respect the fastest transport even if the corresponding transported component is not present in the domain. An illustrative representative numerical example is given later in this paper.

As the governing partial differential equations belong formally to the class of nonlinear hyperbolic systems, it is well-known that the solution of such problems can develop discontinuities (shocks) even for initially smooth profiles. Moreover, numerical methods in the presence of such non-smooth solutions will develop unphysical oscillations that can completely destroy the plausibility of numerical results. Therefore, we extend the second order accurate form of the scheme to the high-resolution form that can locally adjust the stencil of the numerical scheme analogously to (Weighted) Essentially Non-Oscillatory type approximations \cite{shu_essentially_1998,donat2018implicit,zakovafrolkovivc2026esaim} or limit the approximation locally towards a first order accurate scheme that is necessary only in some situations when time steps is violating the CFL condition \cite{duraisamy_implicit_2007, puppo_quinpi_2022,frolkovic2023high, zakovafrolkovivc2026esaim}.

The paper is organized as follows. Section~\ref{sec-mathmodel} introduces the mathematical model, first for a scalar transport equation with Langmuir sorption and then for the corresponding multicomponent system. Section~\ref{sec-numerical-method_secondorder} presents the first and second order compact implicit schemes together with the analogous explicit discretizations. The high-resolution formulation, including the WENO reconstruction, time limiter, and characteristic Courant numbers, is described in Section~\ref{sec-high-resolution-scheme}. Numerical experiments are presented in Section~\ref{sec-experiments}, including a smooth three-component Langmuir problem, a linear benchmark problem, and a three-component displacement chromatography example. Finally, the conclusion is in the last section.

\section{Mathematical model}
\label{sec-mathmodel}

\subsection{Scalar transport equation with sorption}
\label{sec-mathmodel-scalar}

To introduce the topic more clearly, we firstly consider a representative nonlinear transport scalar problem in a porous medium \cite{brusseau1995effect, kacurSolutionContaminantTransport2005, frolkovic_semi-analytical_2006, frolkovicNumericalSimulationContaminant2016}. Let $c=c(x,t)\geq 0$ denote the concentration of the transported substance, where $t\in(0,T)$ and $x\in(x_L,x_R)\subset R$. The transport of the mobile phase is driven by a given constant velocity $u>0$. The governing equation is written in the conservative form
\begin{equation}
    \label{eq_general}
    F(c)_t +  \left( u c \right)_x = 0 \,,
\end{equation}
with the prescribed initial condition 
\begin{equation}
    \label{initialcondition}
    c(x,0)=c^0(x) \,.
\end{equation}
While $c=c(x,t)\geq 0$ denotes the mobile concentration, $F(c)$ describes the total concentration, including both the mobile and adsorbed phases, as given by
\begin{equation}
    \label{total_concentration}
    F(c) = c + \frac{1-\varepsilon}{\varepsilon}\Psi(c) \,,
\end{equation}
where $\varepsilon\in(0,1)$ denotes the porosity of the medium, i.e. the ratio of the void space that can be filled by the liquid solution to the total volume of the porous medium, and $\Psi(c)$ denotes the sorption isotherm. The factor
\begin{equation}
    \label{porosity_factor}
    R = \frac{1-\varepsilon}{\varepsilon}
\end{equation}
measures the relative contribution of the adsorbed phase to the mobile phase. 

The mathematical model can describe not only the transport of contaminants \cite{brusseau1995effect, kacurSolutionContaminantTransport2005, frolkovic_semi-analytical_2006, frolkovicNumericalSimulationContaminant2016} but also of liquids in chromatography \cite{guMathematicalModelingScaleUp2015, schmidt-traubPreparativeChromatography2020, javeed2011efficient,donat2018implicit,burger_linearly_2018,donat2024weno} with the sorption isotherm. In this work, we consider the nonlinear Langmuir sorption isotherm \cite{guMathematicalModelingScaleUp2015,
schmidt-traubPreparativeChromatography2020,
brusseau1995effect,donat2018implicit, burger_linearly_2018,donat2024weno}
\begin{equation}
    \label{Langmuir}
    \Psi(c) = \frac{a c}{1+b c} \,,     \qquad a,b>0 \,,
\end{equation}
where $a$ characterizes the amount of adsorbed substance and $b$ controls how rapidly the adsorption approaches its saturation value.
Consequently, $F(c)$ is
\begin{equation}
    \label{F_langmuir}
    F(c) = c + \frac{1-\varepsilon}{\varepsilon}  \frac{a c}{1+b c} \,,
\end{equation}
with the derivative of the function
\begin{equation}
    \label{F_derivative}
        F'(c) = 1 + \frac{1-\varepsilon}{\varepsilon} \frac{a}{(1+b c)^2} \,.
\end{equation}
Since $a,b>0$ and $0<\varepsilon<1$, we have
\begin{equation}
    F'(c)>1 \qquad \text{for all } c\geq 0 \,.
\end{equation}

The equation \eqref{eq_general} can also be written in form (if all involved derivatives exist)
\begin{equation}
    F'(c)c_t + u c_x = 0 \,.
\end{equation}
Equivalently, in the primitive variables formulation, or the non-conservative formulation \cite{leveque_finite_2004,frolkovic_semi-analytical_2006}, as
\begin{equation}
    \label{transport_speed_scalar}
    c_t + \frac{u}{F'(c)} c_x = 0 \,,
\end{equation}
where
\begin{equation}
    \label{lambda_scalar}
    \lambda(c)=\frac{u}{F'(c)} \,,
\end{equation}
is the nonlinear characteristic speed.

Since $u>0$ and $F'(c)>1$ for $c\geq 0$, the characteristic speed is positive. Therefore, 
a Dirichlet boundary condition is prescribed only at the inflow boundary $x=x_L$ \cite{leveque_finite_2004},
\begin{equation}
    \label{boundarycondition}
    c(x_L,t)=c_{x_L}(t) \,,
\end{equation}
while no boundary condition is required at the outflow boundary $x=x_R$.

The sorption term introduces nonlinear retardation into the transport process. Larger values of $F'(c)$ correspond to stronger retardation and therefore to smaller characteristic speeds. In particular, increasing the adsorption parameter $a$ or decreasing the porosity $\varepsilon$ increases the retardation effect. On the other hand, if the sorption is weak, then $F'(c)$ is close to one and the characteristic speed is close to the transport velocity $u$.

\subsection{Extension to a system of equations}
\label{sec-mathmodel-system}

We now extend the scalar model to a system of transported concentrations. Let $N_c$ denote the number of transported components. We introduce the vector
\begin{equation}
    \label{concentration_vector}
    C = (c_1,c_2,\ldots,c_{N_c})
\end{equation}
of mobile concentrations, where $c_m=c_m(x,t)\geq 0$, $m = 1,\ldots, N_c$. 
The transport of all components is driven by the same constant velocity $u>0$. With
\begin{equation}
    \label{total_concentration_vector}
    F(C)=\bigl(F_1(C),F_2(C),\ldots,F_{N_c}(C)\bigr) \,,
\end{equation}
where $F_m=F_m(C)$ denotes the total concentration of the $m$-th component, $m = 1,\ldots, N_c$, the governing system is written in the form
\begin{equation}
    \label{eq_system_general}
    F_m(C)_t + ( u c_m )_x = 0,  \qquad m = 1,\ldots, N_c \,,
\end{equation}
or, equivalently, in vector form,
\begin{equation}
    \label{eq_system_vector}
    F(C)_t + (uC)_x=0\,.
\end{equation}

The initial condition is prescribed by
\begin{equation}
    \label{initialcondition_system}
    C(x,0)=C^0(x) \,,
\end{equation}
or equivalently,
\begin{equation}
    \label{initialcondition_system_components}
    c_m(x,0)=c_m^0(x),
    \qquad m = 1,\ldots, N_c\,.
\end{equation}

For the system, we again consider a Langmuir sorption isotherm, with $F_m(C)$ defined as
\begin{equation}
    \label{Fi_system}
    F_m(C) = c_m + \frac{1-\varepsilon}{\varepsilon} \frac{a_m c_m}{1+\sum_{j=1}^{N_c} b_j c_j}\,,   \qquad m=1,\ldots, N_c \,.
\end{equation}
This multicomponent form of the Langmuir isotherm is commonly used in mathematical models of nonlinear chromatography \cite{guMathematicalModelingScaleUp2015, javeed2011efficient,
donat2018implicit, burger_linearly_2018, donat2024weno, guiochon2006fundamentals}.

To derive the quasilinear form of the system, we apply the chain rule. Let
\begin{equation}
    \label{Jacobian_definition}
    J(C)=\frac{\partial F}{\partial C}
\end{equation}
be the Jacobian matrix of the mapping $F(C)$, with entries
\begin{equation}
    \label{Jacobian_entries_general}
    J_{mj}(C)=\frac{\partial F_m}{\partial c_j}(C)\,,
    \qquad m,j=1,\ldots, N_c\,.
\end{equation}
Then
\begin{equation}
    \label{chain_rule_system}
    F(C)_t = J(C) C_t \,.
\end{equation}
Therefore, the system \eqref{eq_system_vector} can be written as
\begin{equation}
    \label{system_quasilinear}
    J(C)C_t+uC_x=0 \,.
\end{equation}
If $J(C)$ is invertible, this is equivalent to
\begin{equation}
    \label{system_transport_matrix}
    C_t + uJ(C)^{-1}C_x = 0 \,.
\end{equation}

The matrix $J(C)$ plays the role of a nonlinear retardation matrix. The characteristic speeds of the system are given by the eigenvalues of $J(C)^{-1}$  \cite{donat2018implicit, burger_linearly_2018, donat2024weno,mazzotti2009nonclassical} multiplied by velocity $u$.
If $\mu_m(C)$, $m = 1,\ldots, N_c$, are the eigenvalues of $J(C)$, then the corresponding characteristic speeds are
\begin{equation}
    \label{system_speeds}
    \lambda_m(C)=\frac{u}{\mu_m(C)},
    \qquad m = 1,\ldots, N_c\,.
\end{equation}
Thus, large eigenvalues of $J(C)$ correspond to strong retardation and slow wave propagation, while eigenvalues close to one correspond to weak retardation and faster transport \cite{donat2018implicit,burger_linearly_2018}.

In the Langmuir case considered here, the eigenvalues of $J(C)$, and also $J(C)^{-1}$, are positive \cite{donat2018implicit, donat2024weno, burger_linearly_2018, donat2010secular}, and therefore the direction of propagation is determined by the sign of the transport velocity $u$. In all numerical experiments presented here, we assume $u>0$, so that all characteristic waves propagate from left to right. 

Since the velocity and eigenvalues are all positive, boundary conditions are prescribed only at the inflow boundary $x=x_L$,
\begin{equation}
    \label{boundarycondition_system}
    C(x_L,t)=C_{x_L}(t) \,,
\end{equation}
while no boundary condition is imposed at the outflow boundary $x=x_R$.

\section{Second order accurate schemes}
\label{sec-numerical-method_secondorder}

In this section, we briefly describe the numerical schemes used in the computations. Our aim is to highlight the advantages of several variants of the compact implicit scheme when applied to this kind of problems. The detailed derivation of the compact implicit high-resolution method was presented in \cite{frolkovic2023high, zakovafrolkovic2025tatra}. Here, we recall only the main information and adapt the process to the multicomponent Langmuir system considered in this paper.

The spatial discretization is based on the finite volume method \cite{leveque_finite_2004,harten1997high}. The computational domain is divided into finite volumes $I_i=(x_{i-1/2},x_{i+1/2})$, $i=1,\ldots,M$, with uniform mesh size $h$. The time levels are denoted by $t^n=n\tau$, where $\tau$ is the time step and $n = 0, \ldots, N$. For the system \eqref{eq_system_general}, we denote by $C_i^n$ the numerical approximation of the vector of concentrations associated with the finite volume $I_i$ at time $t^n$. The corresponding accumulated concentration is denoted by $Q_i^n=F(C_i^n)$.

A finite volume discretization of the governing system over $I_i$ and the time interval $(t^n,t^{n+1})$ leads to
\begin{equation}
    \label{finite_volume_scheme_system}
      Q_i^{n+1} - Q_i^n + \frac{\tau}{h} u  \left( C_{i+1/2}^{n+1/2} - C_{i-1/2}^{n+1/2} \right) = 0 \,,
\end{equation}
where $C_{i+1/2}^{n+1/2}$ denotes a suitable approximation of the time-averaged concentration at the cell interface. Since in this work we assume $u>0$, all numerical fluxes are evaluated using left-biased upwind stencils \cite{leveque_finite_2004, harten1997high}. For cases when $u<0$, see \cite{zakovafrolkovic2025tatra}.

\subsection{Compact implicit scheme}
\label{subsec-second_order_compact_implicit}

For comparison purposes, we use the first order implicit upwind scheme. In this case, the interface value is approximated simply by $C_{i+1/2}^{n+1/2}=C_i^{n+1}$.
The resulting scheme is robust but introduces significant numerical diffusion, especially near sharp concentration fronts \cite{leveque_finite_2004,harten1997high}, see numerical experiments later.

The second order compact implicit scheme is obtained by replacing the interface value with a compact second order approximation involving values from the new and previous time levels \cite{zakovafrolkovivc2026esaim, zakovafrolkovic2025tatra, frolkovic2023high}. In the scalar case, the detailed derivation follows from a Taylor expansion in space and time. In the present system case, the same construction is applied componentwise to the vector $C$
\begin{equation}
    \label{C_approx_2ndord_implicit}
        C_{i+1/2}^{n+1/2} = C_{i}^{n+1} - \frac{1}{2}\big(\omega(C_{i-1}^{n+1} - C_{i}^{n}) + (1-\omega)(C_{i}^{n+1} - C_{i+1}^{n})\big) \,,
\end{equation}
with $\omega \in [0,1]$. 

It is important to note that the resulting compact implicit scheme given by \eqref{finite_volume_scheme_system} with \eqref{C_approx_2ndord_implicit} contains only unknowns $C_i^{n+1}$, $C_{i-1}^{n+1}$, and $C_{i-2}^{n+1}$ that we use to apply efficient algebraic solvers, see later. 

For nonsmooth solutions, the second order approximation may produce oscillations near discontinuities \cite{leveque_finite_2004, shu_essentially_1998,duraisamy_implicit_2007}. Therefore, we also use a high-resolution version of the method \cite{frolkovic2023high,  zakovafrolkovivc2026esaim, zakovafrolkovic2025tatra}, in which $\omega$ changes for each concentration and finite volume $I_i$, combined with a time limiter, and the scheme then transforms to the form
\begin{equation}
    \label{C_approx_HR_implicit}
    C_{i+1/2}^{n+1/2} = C_i^{n+1} - \frac{1}{2} l_i \circ
    \left[
        \omega_i \circ \left(C_{i-1}^{n+1}-C_i^n\right) + \left(1-\omega_i\right) \circ \left(C_i^{n+1}-C_{i+1}^n\right)
    \right] \,.
\end{equation}
Here, both $\omega_i=(\omega_{1,i},\ldots,\omega_{N_c,i})$ and
$l_i=(l_{1,i},\ldots,l_{N_c,i})$ are vectors whose components are
defined separately for each concentration in the finite volume $I_i$.
The symbol $\circ$ denotes the Hadamard product, and all operations involving $\omega_i$ and $l_i$ are understood componentwise, in particular, $1-\omega_i$ denotes the vector with components $1-\omega_{m,i}$. Parameter vector $\omega$ is then computed using a WENO reconstruction \cite{shu_essentially_1998, shu1988efficient, jiang1996efficient}, see later in Section~\ref{subsec-high-resolution-scheme_compact_implicit}.

In this formulation, both, the WENO approximation and the time limiter $l$ \cite{sweby1984high} are considered as parts of the numerical scheme, but they play different roles. The WENO approximation is used to determine the reconstruction parameter $\omega$ by assigning larger weights to locally smooth stencils and smaller weights to stencils that cross steep gradients or discontinuities. Thus, WENO does not reduce the scheme directly to the first order approximation, instead, it adapts the reconstruction towards a smoother local stencil. The time limiter $l$, on the other hand, controls the size of the entire second order correction. For small Courant numbers or in smooth regions, the time limiter is inactive and should satisfy $l\approx 1$, so that the full second order approximation is retained. For sufficiently large Courant numbers, particularly near steep fronts or discontinuities, the limiter may decrease towards zero and the scheme is locally reduced towards the first order upwind method \cite{duraisamy_implicit_2007, frolkovic2023high, zakovafrolkovivc2026esaim}.

The obtained system of nonlinear algebraic equations is solved iteratively using the fast sweeping method \cite{lozano2021implicit}, where each iteration is given by a nonlinear Gauss-Seidel iteration \cite{zhao2005fast,frolkovic2023high,
zakovafrolkovivc2026esaim}. At each sweep, a local nonlinear algebraic system is solved for every finite volume by Newton's method. Since in the present paper we consider only positive velocities, the Gauss-Seidel iteration is applied sequentially only from the left to the right, i.e., from the inflow boundary to the outflow boundary.

\subsection{Explicit schemes}
\label{subsec-explicit-schemes}

For comparison with the compact implicit schemes, we also implement analogous explicit schemes of first and second order of accuracy. The main motivation for including these methods is to demonstrate that, for the present class of problems, the use of an explicit time discretization is not necessarily computationally advantageous. 

The corresponding finite volume update can be written in the form, for $Q_i^{n+1} = F(C_i^{n+1})$, as
\begin{equation}
    \label{second_order_system_update}
        Q_i^{n+1} = Q_i^n - \frac{\tau}{h}u \left( C_{i+1/2}^{n+1/2} - C_{i-1/2}^{n+1/2} \right) \,,
\end{equation}
where the interface values $C_{i+1/2}^{n+1/2}$ and $C_{i-1/2}^{n+1/2}$ are reconstructed explicitly from the previous time level.

Since the velocity is assumed to be positive, $u>0$, the first order explicit scheme is written in the upwind form
\begin{equation}
    \label{first_order_explicit_system}
     Q_i^{n+1} = Q_i^n - \frac{\tau}{h}u \left( C_i^n-C_{i-1}^n \right) \,.
\end{equation}
The right-hand side is known from the previous time level. However, after the update of $Q_i^{n+1}$, the new concentration vector $C_i^{n+1}$ must still be recovered from the nonlinear relation $F(C_i^{n+1})=Q_i^{n+1}$, which is solved locally by Newton's method in each finite volume. Note that such numerical approach is necessary, as an inverse function $F^{-1}$ is not available in general.

The second order explicit scheme combines a finite volume spatial discretization with a Lax-Wendroff-type time discretization. The latter is based on the Lax-Wendroff-type approach proposed by Qiu and Shu \cite{qiu_finite_2003}. Their original approach is formulated for hyperbolic systems in the standard conservative form. In the present work, we adopt the Lax-Wendroff-type time discretization idea for the system \eqref{eq_system_vector}, where the nonlinear function $F(C)$ appears under the time derivative. Consequently, the replacement of time derivatives by spatial derivatives involves the Jacobian $J(C)=\partial F(C)/\partial C$ and its inverse. The spatial discretization is then adapted to the present finite volume formulation. WENO reconstructions in finite volume formulations are well established \cite{qiu2002construction, shu2009high} and are used here for the reconstruction of the interface values. 

The main idea of the time discretization is to use a Taylor expansion in time and to replace the resulting time derivatives with spatial derivatives using the governing equation.
For the present system, we have $F(C)_t=-uC_x$. Moreover, using the quasilinear form $C_t+uJ(C)^{-1}C_x=0,$ we obtain $F(C)_{tt}=-u(C_t)_x=u^2\left(J(C)^{-1}C_x\right)_x.$

Using these relations, the second order Taylor approximation in time,
$Q_i^{n+1}\approx Q_i^n+\tau (Q_t)_i^n+\frac{\tau^2}{2}(Q_{tt})_i^n$,
leads to
\begin{align}
    \label{explicit_second_order_taylor}
    Q_i^{n+1}  \approx  \,  &Q_i^n - \tau u C_x \\ 
    \label{explicit_second_order_taylor_secondord}
    & + \frac{\tau^2u^2}{2}
    \left(J(C)^{-1}C_x\right)_x \,.
\end{align} And the corresponding second order explicit scheme is written in the form
\begin{equation}
    \label{second_order_explicit_system_final}
    Q_i^{n+1} = Q_i^n - \tau u D_i^n + \frac{\tau^2u^2}{2}B_i^n \,,
\end{equation}
where $D_i^n$ and $B_i^n$ approximate the first and second spatial terms in \eqref{explicit_second_order_taylor}-\eqref{explicit_second_order_taylor_secondord}, respectively. Their construction is described below using reconstructed values from the previous time level. 

For the vector $C$, we define
\begin{equation}
    \label{explicit_G_reconstruction_component}
    C_{i+1/2}^n = C_{i}^n - \frac{1}{2}
    \left[ \omega \left(C_{i-1}^n - C_{i}^n\right) +
    (1-\omega) \left(C_{i}^n - C_{i+1}^n\right) \right] \,,
\end{equation}
with $\omega\in[0,1]$. The first spatial derivative in
\eqref{explicit_second_order_taylor} is then approximated as 
\begin{equation}
    \label{explicit_first_derivative_G}
    D_i^n = \frac{C_{i+1/2}^n - C_{i-1/2}^n}{h} \,.
\end{equation}

The spatial term $\left(J(C)^{-1}C_x\right)_x$ entering the second order correction \eqref{explicit_second_order_taylor_secondord} is approximated by
$B_i^n$, defined as
\begin{equation}
    \label{explicit_B_definition}
    B_i^n = \frac{A_{i+1/2}^{n}-A_{i-1/2}^{n}}{h} \,,
\end{equation}
where $A_{i+1/2}^{n}$ and $A_{i-1/2}^{n}$ approximate $J(C)^{-1}C_x$ at the corresponding locations
\begin{align}
    \label{explicit_A_plus_minus}
    A_{i+1/2}^{n} &= J(C_{i+1/2}^{n})^{-1}D_{i+1/2}^{n} \,,
    &
    A_{i-1/2}^{n} &= J(C_{i-1/2}^{n})^{-1}D_{i-1/2}^{n} \,.
\end{align}
The corresponding approximations of $C_x$ are
\begin{align}
    \label{explicit_D_plus_minus}
    D_{i+1/2}^{n} &= \frac{C_{i+3/2}^{n}-C_{i-1/2}^{n}}{2h} \,,
    &
    D_{i-1/2}^{n} &= \frac{C_{i+1/2}^{n}-C_{i-3/2}^{n}}{2h} \,.
\end{align}

The concentration vector $C_i^{n+1}$ is then recovered from $Q_i^{n+1} = F(C_i^{n+1})$ in \eqref{second_order_explicit_system_final},  which is solved locally by Newton's method in each finite volume.

Compared with the first order explicit scheme, the second order explicit scheme uses a wider stencil because of the reconstruction of interface values and the approximation of spatial derivatives. Therefore, additional neighboring values, or ghost values near the boundary, are required, introduced by extrapolation so that the full stencil can be evaluated. Moreover, despite being explicit with respect to the transport term, the scheme still requires a local Newton iteration due to the nonlinear relation between $F(C)$ and $C$. In addition, it remains subject to the usual stability restriction on the time step. This makes it unsuitable for the large Courant numbers considered in the numerical experiments below, where the compact implicit schemes can be used with much larger time steps \cite{leveque_finite_2004, qiu_finite_2003, duraisamy_implicit_2007}.

\section{High-resolution schemes}
\label{sec-high-resolution-scheme}

\subsection{High-resolution compact implicit scheme}
\label{subsec-high-resolution-scheme_compact_implicit}

We now continue with the compact implicit numerical scheme in its high-resolution form. For constant $u>0$, the high-resolution scheme \cite{frolkovic2023high,zakovafrolkovivc2026esaim} in vector notation adapted to our problem is written using \eqref{C_approx_HR_implicit} as
\begin{align}
    \label{HR_implicit_fullscheme}
        Q_i^{n+1} - Q_i^n +     
       \frac{\tau}{h}u \bigg[ \Big( &C_{i}^{n+1}   - \frac{1}{2} l_i \circ \big(\omega_{i} \circ   (C_{i-1}^{n+1} - C_{i}^{n}  ) + (1-\omega_i) \circ   (C_{i}^{n+1}   - C_{i+1}^{n}) \Big) \nonumber \\
     -  \Big(&C_{i-1}^{n+1} - \frac{1}{2} l_{i-1} \circ\big(\omega_{i-1} \circ (C_{i-2}^{n+1} - C_{i-1}^{n}) + (1-\omega_{i-1}) \circ (C_{i-1}^{n+1} - C_{i}^{n}  ) \Big) \bigg]  = 0 \,.
\end{align}
As mentioned in Section \ref{subsec-second_order_compact_implicit} $\omega_i=(\omega_{1,i},\ldots,\omega_{N_c,i})$ and $l_i=(l_{1,i},\ldots,l_{N_c,i})$ are vectors.
The values $\omega_{m,i}$ are obtained from the WENO approximation \cite{shu1988efficient, jiang1996efficient, shu_essentially_1998}. A detailed derivation of the WENO approximation for the compact implicit scheme has already been presented in \cite{frolkovic2023high, zakovafrolkovivc2026esaim, zakovafrolkovic2025tatra} and is therefore not repeated here, while the values $l_{m,i}$ are computed from the limiter formula \cite{frolkovic2023high} for each component $m$ as
\begin{equation}
    \label{limiter} 
    l_{m,i} = \min \Biggl{\{} 1 , \max \Biggl[ 0, \left( \omega_{m,i} + \frac{1-\omega_{m,i} }{r_{m,i}}\right)^{-1}\left(\frac{2}{ \mathrm{Co}_{m,i} }+l_{m,i-1}(\omega_{m,i-1} r_{m,i-1} + 1-\omega_{m,i-1})\right)\Biggr] \Biggr{\}}\,, \quad m = 1, \ldots , N_c \,,
\end{equation}
with 
\begin{equation}
    \label{ratio_r}
    r_{m,i}=\frac{ c_{m,i-1}^{n+1} - c_{m,i}^{n} }{c_{m,i}^{n+1} - c_{m,i+1}^{n}} \,.
\end{equation}
The time limiter is important only for computations with large time steps that are typically controlled. It controls the size of the second order correction in \eqref{HR_implicit_fullscheme}. In the scalar case, the limiter is constructed using a local Courant number ($\mathrm{Co}_{m,i}$) (defined later in text as $\mathrm{Co}_{m,i}^{\mathrm{char}}$), which enters the condition used to preserve the non-oscillatory behavior of the high-resolution scheme, see \cite{frolkovic2023high, zakovafrolkovivc2026esaim}. In the present vector case, the choice of this Courant number requires additional attention.

For the scalar Langmuir isotherm \eqref{F_langmuir} we have the derivative \eqref{F_derivative}, for which the retardation factor is always greater than one. Consequently, the characteristic speed of the scalar equation is reduced according to $\lambda(c)=\frac{u}{F'(c)}$. Since $F'(c)>1$, we have $\lambda(c) < u$. Therefore, we can define
\begin{equation}
    \label{Courantnumber_max}
        \mathrm{Co}_{\max} = \frac{\tau}{h}u \,,
\end{equation}
which gives an upper bound for scalar Courant number.

For systems, the role of the retardation factor is played by the Jacobian matrix $J(C)=\partial F(C)/\partial C$. The characteristic speeds are given by the eigenvalues of $J(C)^{-1}$, multiplied by velocity $u$. 

In all numerical experiments presented here, we assume $u>0$, so that all characteristic waves propagate from left to right. The corresponding numerical stencils are therefore chosen in the upwind direction. The case of negative transport velocity for scalar equations \cite{zakovafrolkovivc2026esaim} with Freundlich isotherms, including the corresponding modification of the computational algorithm, was discussed in \cite{zakovafrolkovic2025tatra}.

Since the eigenvalues of $J(C)$ are positive, and considering the Langmuir sorption, they are larger than or equal to $1$. The eigenvalues of $J(C)^{-1}$ are also positive and fall within the range of $(0,1]$. A large eigenvalue of $J(C)$ corresponds to stronger retardation and hence to a smaller characteristic speed, whereas a small eigenvalue of $J(C)$ corresponds to weaker retardation and faster transport. Therefore, the maximum Courant number can be bounded, as in \eqref{Courantnumber_max}, or in terms of the largest eigenvalue of $J(C)^{-1}$. In the numerical experiments below, we first use such an upper bound to choose the time step. Later, for selected test cases, we compute the Courant number more precisely from the eigenvalues of $J(C)^{-1}$, used for the computation of limiter \eqref{limiter}.

If the upper bound $\mathrm{Co}_{\max}$ is used directly in the limiter formula for system of equations, then the limiting may be too strong for some components. In that case, the values $l_{m,i}$ may be reduced unnecessarily, and the high-resolution scheme is driven too close to the first order upwind method. This increases numerical diffusion and may reduce the accuracy of the solution, especially in problems where the actual characteristic speeds are significantly smaller than the velocity $u$ because of sorption.

For this reason, in selected computations we use characteristic Courant numbers based on the eigenvalues of $J(C)^{-1}$. Let $\lambda_{m,i}$, $m=1,\ldots,N_c$, denote the eigenvalues of $J(C_i)^{-1}$ in the finite volume $I_i$. The corresponding characteristic Courant numbers are defined by
\begin{equation}
    \label{Courantnumber_characteristic_components}
    \mathrm{Co}_{m,i}^{\mathrm{char}} = \frac{\tau u}{h}\left|\lambda_{m,i}\right| \,,
    \qquad m=1,\ldots,N_c \,.
\end{equation}
The maximum characteristic Courant number associated with the finite volume $I_i$ is then
\begin{equation}
    \label{Courantnumber_characteristic_local}
    \mathrm{Co}_{i,\max}^{\mathrm{char}} = \max_{m=1,\ldots,N_c} \mathrm{Co}_{m,i}^{\mathrm{char}}
    = \frac{\tau u}{h} \rho\left(J(C_i)^{-1}\right) \,,
\end{equation}
where $\rho$ denotes the spectral radius. The Courant number used in the limiter can therefore be chosen as the local value $\mathrm{Co}_{i,\max}^{\mathrm{char}}$ instead of the global upper bound $\mathrm{Co}_{\max}$ defined in \eqref{Courantnumber_max}.

This choice gives a less restrictive limiter, because it reflects the actual wave speeds present in the numerical solution. If the sorption strongly retards the transport, then the eigenvalues of $J(C)^{-1}$ are significantly smaller than 1, alternatively, eigenvalues of $uJ(C)^{-1}$ are smaller than $u$, and the corresponding characteristic Courant numbers are smaller than the global bound. Consequently, the limiter does not reduce the second order correction more than necessary.

In the numerical experiments below, we therefore distinguish between two uses of the Courant number. The time step is chosen using the upper bound $\mathrm{Co}_{\max}=\tau u/h$, which gives a simple and safe estimate. However, when the limiter is evaluated, we may use the local characteristic Courant numbers $\mathrm{Co}_{i,\max}^{\mathrm{char}}$ computed from \eqref{Courantnumber_characteristic_local}. This prevents excessive limiting and allows the high-resolution scheme to preserve more of its second order accuracy while still controlling oscillations near steep fronts or discontinuities.

\subsection{WENO approximation in explicit scheme}
\label{subsec-explicit_WENO}

For the explicit scheme, parameter $\omega$ in \eqref{explicit_G_reconstruction_component} may change according to the finite volume and is replaced componentwise by $\omega_{m,i}$ for each concentration $m$ and finite volume $I_i$, and is computed from values at previous time levels. Then $\omega_{m,i}$ can either be fixed or determined adaptively using a WENO procedure. Following again the standard WENO approach
\cite{shu1988efficient, jiang1996efficient, shu_essentially_1998}, also used in \cite{zakovafrolkovic2025tatra, zakovafrolkovivc2026esaim}, we define the upwind and downwind differences for each component by
\begin{equation}
    \label{explicit_WENO_differences}
    \Delta_{m,i}^{\mathrm{u},n} = c_{m,i}^n-c_{m,i-1}^n \,,
    \qquad
    \Delta_{m,i}^{\mathrm{d},n} = c_{m,i+1}^n-c_{m,i}^n \,.
\end{equation}
The corresponding nonlinear weights are
\begin{equation}
    \label{explicit_WENO_weights}
    a_{m,i}^{\mathrm{u},n} = \frac{\bar{\omega}} {\left( \epsilon_w + \left(\Delta_{m,i}^{\mathrm{u},n}\right)^2
    \right)^2}\,,  \qquad 
    a_{m,i}^{\mathrm{d},n} = \frac{1-\bar{\omega}} {\left( \epsilon_w +
        \left(\Delta_{m,i}^{\mathrm{d},n}\right)^2 \right)^2} \,,
\end{equation}
where $\epsilon_w>0$ is a small parameter to avoid division by zero and $\bar{\omega}$ is the preferred value. In the computation, we use $\epsilon_w = 10^{-10}$ and $\bar{\omega}=\frac{1}{2}$. The WENO reconstruction parameter is then obtained by
\begin{equation}
    \label{explicit_WENO_omega}
    \omega_{m,i} = \frac{a_{m,i}^{\mathrm{u},n}} {a_{m,i}^{\mathrm{u},n}+a_{m,i}^{\mathrm{d},n}}\,,  \qquad m=1,\ldots,N_c \,.
\end{equation}
In smooth regions with comparable differences, $\omega_{m,i}$ remains close to $\bar{\omega}=1/2$. 

No additional time limiting is applied in the explicit scheme, since the explicit time discretization is subject to a stability restriction on the Courant number, sufficiently small time steps have to be used, and the time limiter introduced for the compact implicit scheme to control the second order correction at large Courant numbers is therefore not required.

\section{Numerical experiments}
\label{sec-experiments}

In this section, most of the experiments concern multicomponent models used in chromatography  \cite{javeed2011efficient,donat2018implicit,  burger_linearly_2018,donat2024weno}. We present numerical experiments for models with the Langmuir isotherm \eqref{Langmuir}. Different choices of the parameters $a_m$, $b_m$, $(a_m,b_m>0)$ for $m = 1,\ldots, N_c$ are considered in order to test the accuracy and efficiency of the proposed numerical schemes.

In addition, a linear hyperbolic benchmark with a known characteristic decomposition is included \cite{leveque_finite_2004}.

To evaluate the accuracy of the methods, we compute the discrete $L_1$-error ($E$) at the final time $T$. For systems, the error is computed component-wise and summed over all components as
\begin{equation}
    \label{error}
    E = h \sum_{i=1}^{M}\left(\left| c_{1,i}^N - \bar{c}_{1,i}^N \right|+  \ldots + \left| c_{N_c,i}^N - \bar{c}_{N_c,i}^N \right|\right) \,,
\end{equation}
with $\bar c$ denoting the exact or reference solution. 
The experimental order of convergence $(EOC)$ is then calculated from two consecutive mesh refinements for the examples with the known exact solution or reference solution. 

For each of the experiments using Langmuir sorption, we state the upper bound of Courant number $\mathrm{Co}_{\max}$, based on \eqref{Courantnumber_max}. This value is used as a convenient global measure of the time-step size and gives a safe upper bound based only on the transport velocity $u$. When the high-resolution scheme is used and the limiter is active, we also compute the local characteristic Courant numbers, as described in Section~\ref{subsec-high-resolution-scheme_compact_implicit}. These values do not exceed the global upper bound $\mathrm{Co}_{\max}$ in the Langmuir cases considered here, and they are often smaller. For the linear benchmark, the characteristic Courant numbers are considered separately.

The numerical methods were implemented in the Python programming language \cite{CS-R9526}, using the NumPy library \cite{harris2020array} for numerical linear algebra and array operations, and the graphical outputs were produced using Matplotlib \cite{hunter2007matplotlib}.

\subsection{First experiment: smooth three-component Gaussian profiles}
\label{exp1}

In this experiment, inspired by the test problem considered in \cite{donat2018implicit, donat2024weno}, we study a smooth three-component problem. The aim is to verify experimentally that the proposed second order schemes achieve the expected order of accuracy and to compare their computational efficiency. 

We prescribe smooth initial data in the form of three Gaussian profiles,
\begin{equation}
    \label{experiment1_concentration}
    c_m = \varrho_m \exp(-100(x-1/2)^2)\,, \quad m=1,2,3 \,,
\end{equation}
with $\varrho_1=1,\, \varrho_2=2,\, \varrho_3=3.\,$ The parameters in Langmuir isotherm \cite{donat2018implicit,donat2024weno} are set to $a_1=4,\, a_2=5,\, a_3=6,\,b_1=b_2=b_3=1$, the velocity $u=0.2$, the porosity parameter is set to $\varepsilon=0.5$, and we compute the numerical solutions on a sequence of refined meshes.

\begin{figure}[H]
    \begin{center}
        \includegraphics[width=0.8\paperwidth]{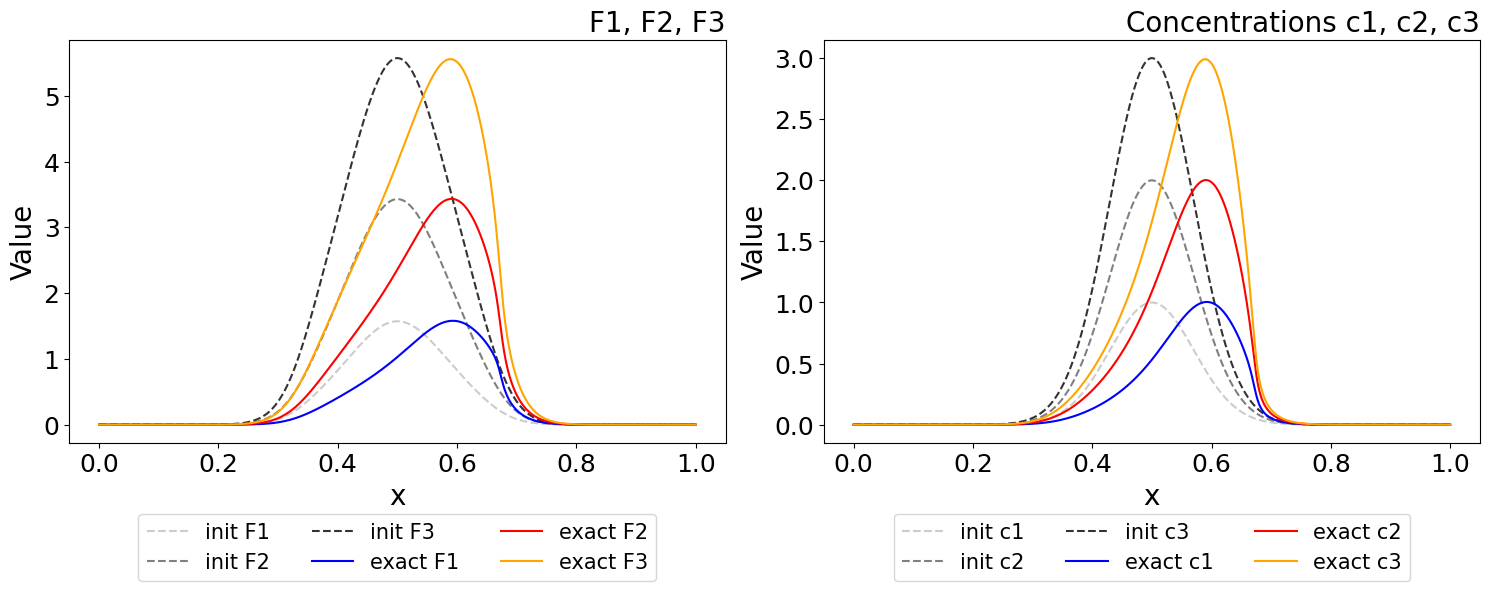}
        \caption{Initial profiles at $t=0$, see eq. \eqref{experiment1_concentration}, and the corresponding reference solution at $T=0.5$ for the smooth three-component Langmuir problem (Experiment \ref{exp1}). The left figure shows the accumulated concentrations $F_1$, $F_2$, and $F_3$, while the right panel shows the mobile concentrations $c_1$, $c_2$, and $c_3$. The reference solution was computed using  the second order compact implicit scheme with the constant reconstruction parameter $\omega=1/2$ and $M_{\mathrm{ref}}=20480$.}
        \label{FIG:exact_gaussians}
    \end{center}
\end{figure}

The reference solution was computed on a fine mesh with $M_{\mathrm{ref}}=20480$ up to the final time $T=0.5$, for $x\in[0,1]$. For this purpose, we used the implicit second order scheme with a constant parameter $\omega$. This reference solution was then used for the computation of the discrete errors and the corresponding experimental orders of convergence. The Figure \ref{FIG:exact_gaussians} shows the comparison between the initial condition and reference solution, for both $c_1, c_2, c_3$ (right figure) and $F_1,F_2,F_3$ (left figure).

The error is measured in the discrete $L_1$-norm \eqref{error} and is used to compute the experimental order of convergence. We compare first order implicit and explicit schemes with second order schemes. For the second order discretization, two variants are considered: one with a constant value of the parameter $\omega$ set to $1/2$, and one using the WENO approximation. Both implicit and explicit time discretizations are tested.

Since the three components propagate with different characteristic speeds, the corresponding Courant numbers are different, but smaller than the upper bound $\mathrm{Co}_{\max}$. In the comparison involving explicit schemes, the stability restriction is imposed by the explicit WENO discretization. For this reason, we use $\mathrm{Co}_{\max} = 0.5,\,$ which is the largest acceptable value in our computations. Moreover, for such a small Courant number, the additional time limiting of the high-resolution schemes is not engaged and is not needed in this test.

The results reported in Table~\ref{TAB:Results_compareEOC_gaussians} show decreasing errors under mesh refinement and confirm the expected convergence rates. 
For the schemes with the constant parameter $\omega=1/2$, the implicit scheme gives smaller errors than the corresponding explicit scheme. For the WENO variants, the errors are comparable, with the explicit scheme giving slightly smaller errors on the finer meshes.

\begin{table}[h!]
    \begin{center}
    \begin{tabular}{c c | c c c | c c c }
    \hline
    \multicolumn{2}{c}{} & \multicolumn{6}{c}{First order}  \\
    \multicolumn{2}{c}{} & \multicolumn{3}{c}{Explicit} & \multicolumn{3}{c}{Implicit}  \\
    \hline
    $M $ & $N$ & $E$                   & EOC  & CPU time  & $E$                    & EOC  & CPU time  \\ 
    \hline
    160  & 32  & 2.6 $\cdot 10^{-2}$ & -    & 0.75   & 5.8 $\cdot 10^{-2}$ & -    & 0.90    \\
    320  & 64  & 1.3 $\cdot 10^{-2}$ & 0.96 & 2.30   & 3.0 $\cdot 10^{-2}$ & 0.92 & 3.07    \\
    640  & 128 & 6.8 $\cdot 10^{-3}$ & 0.98 & 7.19   & 1.5 $\cdot 10^{-2}$ & 0.96 & 11.22   \\
    1280 & 256 & 3.4 $\cdot 10^{-3}$ & 0.99 & 28.22   & 7.9 $\cdot 10^{-3}$ & 0.98 & 37.29  \\
    \hline
    \multicolumn{2}{c}{} &  \multicolumn{6}{c}{Second order with constant $\omega=1/2$}  \\
    \multicolumn{2}{c}{} & \multicolumn{3}{c}{Explicit} & \multicolumn{3}{c}{Implicit}  \\
    \hline
    $M $ & $N$  & $E$                    & EOC  & CPU time & $E$                  & EOC  & CPU time  \\ 
    \hline
    160  & 32  & 1.7 $\cdot 10^{-3}$ & -    & 0.88   & 1.2 $\cdot 10^{-3}$ & -    & 1.82     \\
    320  & 64  & 4.5 $\cdot 10^{-4}$ & 1.89 & 2.97   & 2.9 $\cdot 10^{-4}$ & 2.06 & 4.91     \\
    640  & 128 & 1.2 $\cdot 10^{-4}$ & 1.94 & 11.78  & 7.2 $\cdot 10^{-5}$ & 2.03 & 15.39    \\
    1280 & 256 & 3.0 $\cdot 10^{-5}$ & 1.97 & 43.42  & 1.8 $\cdot 10^{-5}$ & 2.01 & 62.41    \\
    \hline 
    \multicolumn{2}{c}{} &  \multicolumn{6}{c}{Second order with WENO}  \\
    \multicolumn{2}{c}{} & \multicolumn{3}{c}{Explicit} & \multicolumn{3}{c}{Implicit}  \\
    \hline
    $M $ & $N$  & $E$                    & EOC  & CPU time & $E$                  & EOC  & CPU time  \\ 
    \hline
    160  & 32  & 4.9 $\cdot 10^{-3}$ & -    & 1.31   & 4.8 $\cdot 10^{-3}$ & -    & 2.50     \\
    320  & 64  & 1.2 $\cdot 10^{-3}$ & 1.99 & 3.96   & 1.3 $\cdot 10^{-3}$ & 1.82 & 5.26     \\
    640  & 128 & 3.2 $\cdot 10^{-4}$ & 1.93 & 12.08  & 3.8 $\cdot 10^{-4}$ & 1.85 & 22.39    \\
    1280 & 256 & 7.9 $\cdot 10^{-5}$ & 2.00 & 46.57  & 9.5 $\cdot 10^{-5}$ & 1.99 & 92.71    \\
    \hline 
    \end{tabular}
    \end{center}
    \caption{A comparison of the numerical errors $E$, experimental orders of convergence (EOC), and mean CPU times in seconds for the explicit and implicit first order schemes, the second order schemes with the constant reconstruction parameter $\omega=1/2$, and the second order schemes with WENO reconstruction in Experiment~\ref{exp1}. 
    The results are reported for the smooth three-component Gaussian initial data with $\mathrm{Co}_{\max}=0.5$, where $M$ denotes the number of finite volumes and $N$ the number of time steps. The CPU times are arithmetic means obtained from five runs.}
    \label{TAB:Results_compareEOC_gaussians}
\end{table}

\begin{figure}[h]
    \begin{center}
        \includegraphics[width=0.7\paperwidth]{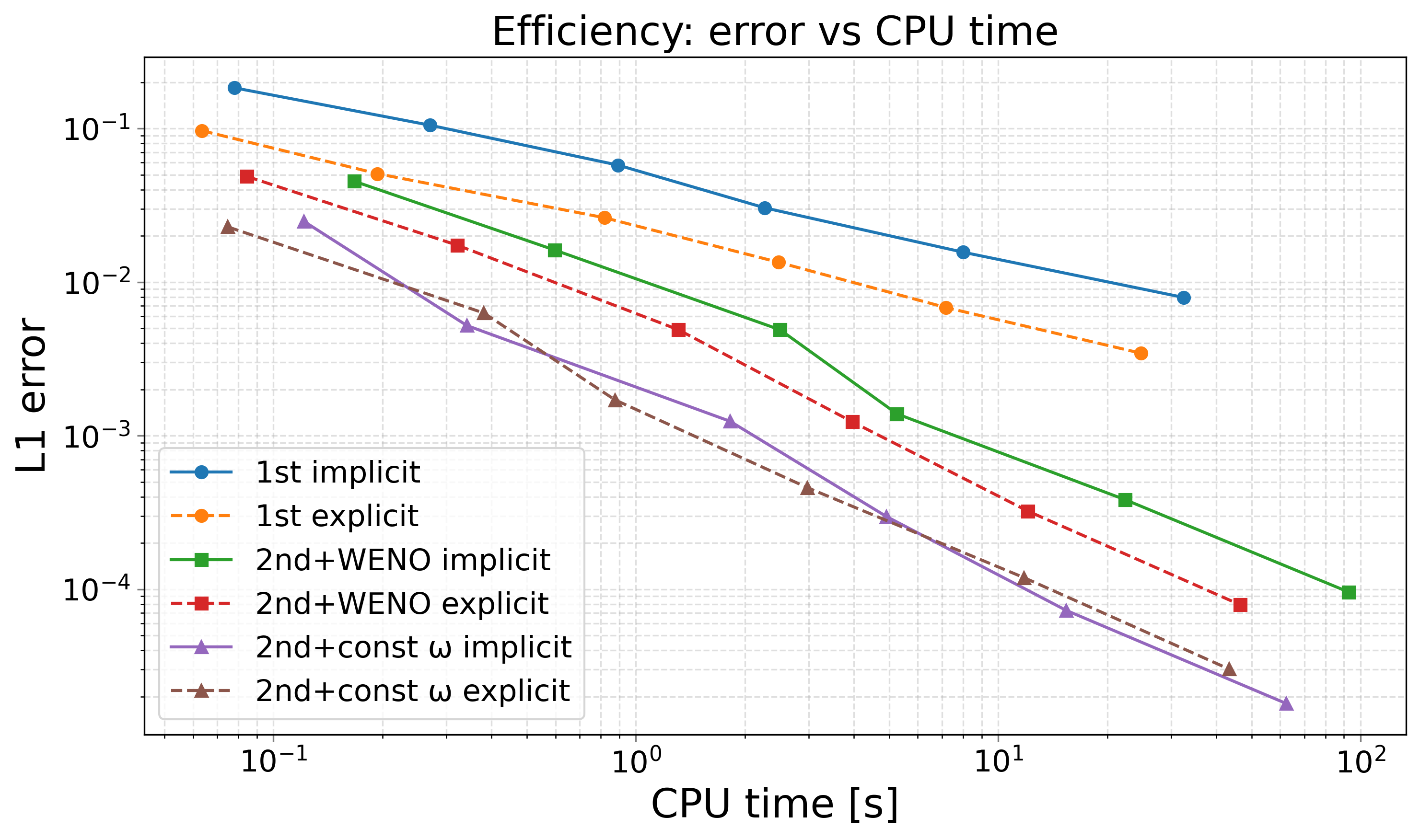}
        \caption{Efficiency comparison of the numerical schemes for Experiment~\ref{exp1} in terms of the $L_1$ error as a function of CPU time. The results are shown on a $log-log$ scale for the explicit and implicit first order schemes, the second order schemes with WENO reconstruction, and the second order schemes with the constant reconstruction parameter $\omega=1/2$. The computations were performed with $\mathrm{Co}_{\max}=0.5$, and the CPU times represent arithmetic means over five runs.}
        \label{FIG:logerror_gaussians}
    \end{center}
\end{figure}

The CPU times provide an additional comparison of the methods. As expected, implicit schemes are more computationally expensive, especially when combined with the WENO reconstruction. The additional computational cost of the implicit formulation is considerably more pronounced for the WENO variant than for the second order scheme with constant $\omega=1/2$. Figure~\ref{FIG:logerror_gaussians} presents a work-precision comparison, where the error is plotted as a function of CPU time on a log-log scale. The figure compares first order schemes, second order schemes with WENO reconstruction, and second order schemes with constant $\omega=1/2$, for both implicit and explicit time discretizations. 
The work-precision comparison shows that the relative efficiency depends on the reconstruction. For the schemes with constant $\omega=1/2$, the implicit method provides smaller errors at a moderate additional computational cost, whereas the explicit WENO scheme is more efficient than the corresponding implicit WENO variant for the small Courant number considered in this test.

The reported CPU times in Table~\ref{TAB:Results_compareEOC_gaussians} (visualized in Figure~\ref{FIG:logerror_gaussians}) were obtained by repeating each simulation five times and taking the arithmetic mean of the measured values. This averaging was used to reduce the influence of random fluctuations in the computational environment.

Although the explicit schemes are computationally competitive for the small Courant number considered in this experiment, their time step remains restricted by the stability condition. We therefore focus in the following computations on the implicit discretizations and repeat the convergence study with the larger Courant number $\mathrm{Co}_{\max}=8$ in order to exploit the weaker time-step restriction of the implicit approach. The corresponding errors and experimental orders of convergence are shown in Table~\ref{TAB:Results_compareEOC_gaussians_highCn}.
In contrast to the computations with $\mathrm{Co}_{\max}=0.5$, for the high-resolution WENO variant, the time limiter is active in this test in order to control the second order correction and suppress spurious oscillations for the larger time steps. The presented results confirm that the implicit schemes can provide stable numerical solutions even for large Courant numbers and preserve, for sufficiently fine grids, the expected EOCs.

\begin{table}[h!]
    \begin{center}
    \begin{tabular}{c c | c c c c c c }
    \hline
    \multicolumn{2}{c}{} & \multicolumn{2}{c}{First order} & \multicolumn{2}{c}{Second order, $\omega=1/2$} & \multicolumn{2}{c}{HR with WENO}  \\
    \hline
    $M $ & $N$ & $E$                   & EOC  & $E$                   & EOC  & $E$ & EOC   \\ 
    \hline
    160  & 2  & 3.0 $\cdot 10^{-1}$ & -    & 5.3 $\cdot 10^{-2}$ & -    & 5.3 $\cdot 10^{-2}$ & -    \\
    320  & 4  & 1.3 $\cdot 10^{-1}$ & 0.68 & 1.6 $\cdot 10^{-2}$ & 1.70 & 1.8 $\cdot 10^{-2}$ & 1.52 \\
    640  & 8  & 7.5 $\cdot 10^{-2}$ & 0.79 & 4.4 $\cdot 10^{-3}$ & 1.87 & 5.7 $\cdot 10^{-3}$ & 1.68 \\
    1280 & 16 & 4.0 $\cdot 10^{-2}$ & 0.88 & 1.1 $\cdot 10^{-3}$ & 1.94 & 1.6 $\cdot 10^{-3}$ & 1.83  \\
    \hline
    \end{tabular}
    \end{center}
    \caption{A comparison of the numerical errors $E$ and experimental orders of convergence (EOC) for the implicit first order scheme, the compact implicit second order scheme with the constant reconstruction parameter $\omega=1/2$, and the compact implicit high-resolution scheme with WENO reconstruction in Experiment~\ref{exp1}. The results are reported for the smooth three-component Gaussian initial data with $\mathrm{Co}_{\max}=8$.}
    \label{TAB:Results_compareEOC_gaussians_highCn}
\end{table}

\subsection{Second experiment: Linear benchmark problem}
\label{exp2}

As next numerical experiment, we consider a linear system obtained by choosing the accumulation functions in the form
\begin{align}
   F_1(c_1, c_2) &= 0.505 c_1 - 0.495 c_2 \,, \nonumber \\
   F_2(c_1, c_2) &= 0.505 c_2 - 0.495 c_1 \,. 
\end{align}

The corresponding Jacobian matrix is constant,
\begin{equation}
    J =
    \begin{bmatrix}
        0.505 & -0.495 \\
        -0.495 & 0.505
    \end{bmatrix} \,.
\end{equation}

The eigenvalues of $J$ are $1$ and $0.01$. Therefore, for $u=1$, the system contains two characteristic waves propagating with speeds $1$ and $100$, respectively. The associated eigenvectors are $(1,-1)^T$ and $(1,1)^T$, hence, the solution can be written as a combination of two traveling waves associated with the slow and fast characteristic variable \cite{leveque_finite_2004}.

This example is useful as a benchmark because it has an explicit analytical structure while still exhibiting two characteristic components with significantly different propagation speeds. It therefore provides a suitable test for verifying the implementation, evaluating experimental convergence rates, and checking whether the numerical method correctly captures both slow and/or fast waves.

We consider two variants of the linear benchmark problem. The first one uses a smooth initial condition, while the second one contains a discontinuity. This allows us to distinguish between the behavior of the schemes in smooth regions, where no limiting is expected to be necessary, and near discontinuities, where limiting plays an essential role.

For the smooth case, we prescribe a Gaussian profile associated only with the slow characteristic component, while the fast characteristic component is set to zero. The exact solution is therefore obtained by shifting the initial profile with unit speed:
\begin{equation}
    \label{linear_gaussian_solution}
    c_1(x,t)  = \frac{1}{2}\exp\left(-50(x-1-t)^2\right),  \qquad c_2(x,t) = - \, \frac{1}{2}\exp\left(-50(x-1-t)^2\right) \,.
\end{equation}
Since the solution is smooth, no space or time limiting is applied in this test. We compare the computational efficiency of the first order scheme and the second order scheme with constant $\omega=1/2$, comparing error $E$ as a function of CPU time.

The computations are performed for a sequence of increasing Courant numbers. Since the linear system contains two characteristic speeds, two Courant numbers are present in each computation. We denote them by $\mathrm{Co}_{S}$ and $\mathrm{Co}_{L}$, where $\mathrm{Co}_{L}=100\, \mathrm{Co}_{S}$.

The values used in these two experiments are $\mathrm{Co}_{L}=12.5, 25, 50, 100, 200, 400, 800$, as reported in Table~\ref{TAB:2example_compareerrorwithcn_gaussians_new}. These large Courant numbers ($\mathrm{Co}_L$) are admissible for the implicit schemes but not for the explicit schemes, which are subject to a much stronger stability restriction.

\begin{figure}[H]
    \begin{center}
        \includegraphics[width=0.8\paperwidth]{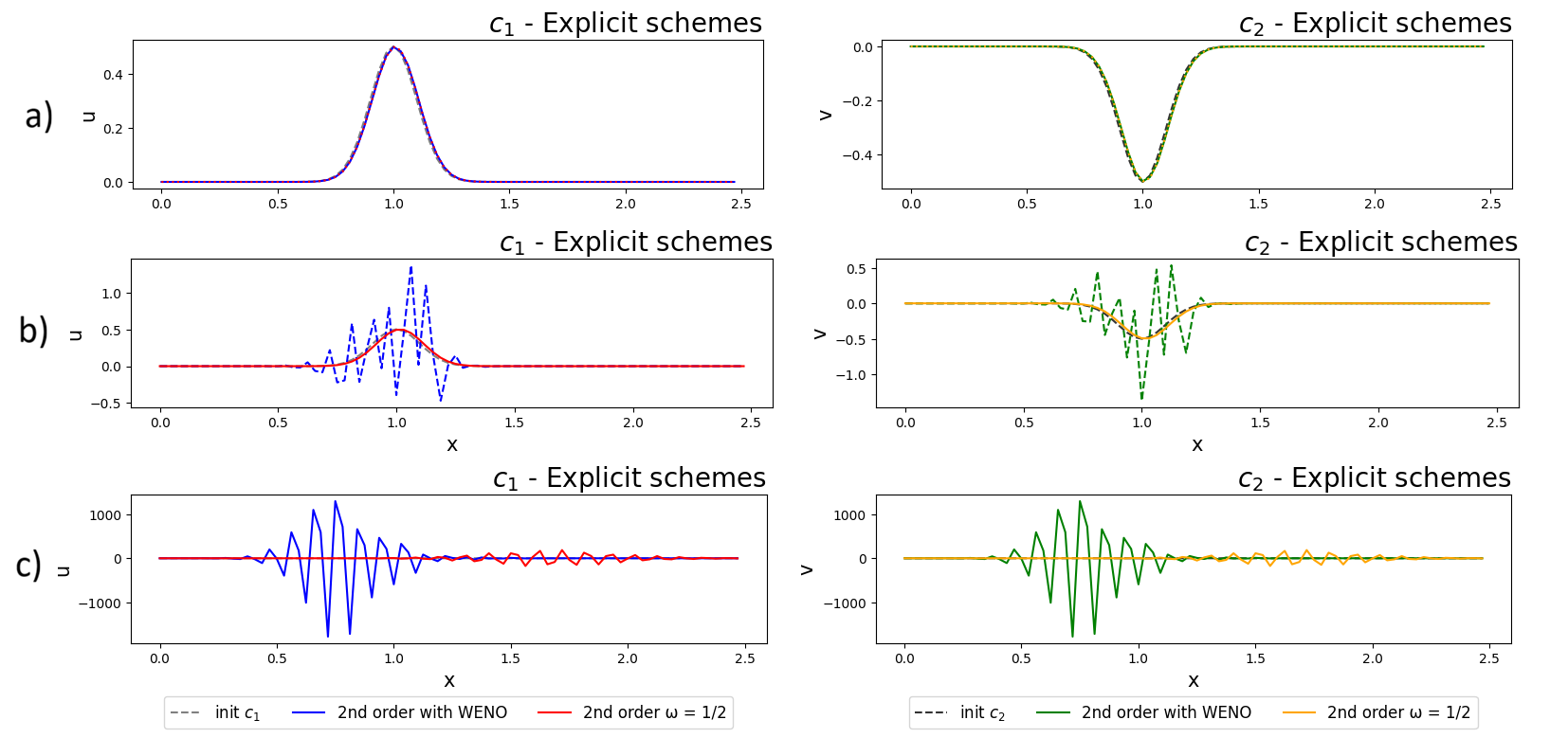}
       \caption{Comparison of two explicit second order schemes, one using the constant reconstruction parameter $\omega=1/2$ and the other using WENO reconstruction, for the smooth case of the linear benchmark problem in Experiment~\ref{exp2}, with the exact Gaussian solution given by \eqref{linear_gaussian_solution}. The results at $t\approx0.15$ are shown for $c_1$ in the left column and for $c_2$ in the right column. From top to bottom, the rows correspond to $\mathrm{Co}_{L}=0.5$, $1$, and $1.25$, denoted by $a)$, $b)$, and $c)$, respectively. Both schemes remain stable for $\mathrm{Co}_{L}=0.5$, oscillations appear in the WENO solution for $\mathrm{Co}_{L}=1$, and both explicit schemes become unstable for $\mathrm{Co}_{L}=1.25$.}
        \label{FIG:gausian_porovnanieexplicitnychschem_rozneCn}
    \end{center}
\end{figure}

This restriction is illustrated in Figure~\ref{FIG:gausian_porovnanieexplicitnychschem_rozneCn}, where explicit second order schemes are compared at approximately $t=0.15$. In the first row, corresponding to $\mathrm{Co}_{L}=0.5$, both the constant $\omega=1/2$ scheme and the scheme with WENO approximation remain stable and non-oscillatory. In the second row, for $\mathrm{Co}_{L}=1$, oscillations already appear in the WENO approximation, while the second order scheme with constant $\omega=1/2$ still remains stable. In the third row, corresponding to $\mathrm{Co}_{L}=1.25$, both explicit second order schemes become unstable. This confirms that explicit schemes are not suitable for the large Courant numbers considered in the efficiency comparison.

Therefore, the following efficiency comparison is restricted to the implicit first order scheme and the compact implicit second order scheme with the constant reconstruction parameter $\omega=1/2$. The corresponding numerical errors $E$ and CPU times for increasing Courant numbers are reported in Table~\ref{TAB:2example_compareerrorwithcn_gaussians_new}.

\begin{table}[h!]    
    \begin{center} 
    \begin{tabular}{c c | c c | c c }
    \hline
    \multicolumn{2}{c}{} & \multicolumn{2}{c}{ } & \multicolumn{2}{c}{Second order compact implicit } \\
    \multicolumn{2}{c}{} & \multicolumn{2}{c}{First order implicit} & \multicolumn{2}{c}{with $\omega=1/2$} \\
    \hline
    $\mathrm{Co}_{L}$ & $\mathrm{Co}_{S}$ &  $E$   & CPU time [s] & $E$   & CPU time [s] \\ 
    \hline
    12.5  & 0.125 & 4.40 $\cdot 10^{-2}$ & 11.16 & 9.58 $\cdot 10^{-4}$ & 12.73   \\
    25    & 0.25  & 4.80 $\cdot 10^{-2}$ & 6.10  & 9.17 $\cdot 10^{-4}$ & 6.51    \\
    50    & 0.50  & 5.56 $\cdot 10^{-2}$ & 3.38  & 7.58 $\cdot 10^{-4}$ & 3.80    \\    
    100   & 1.00  & 6.93 $\cdot 10^{-2}$ & 1.29  & 4.12 $\cdot 10^{-4}$ & 1.51    \\    
    200   & 2.00  & 9.25 $\cdot 10^{-2}$ & 0.60  & 2.94 $\cdot 10^{-3}$ & 0.74    \\
    400   & 4.00  & 1.27 $\cdot 10^{-1}$ & 0.32  & 1.38 $\cdot 10^{-2}$ & 0.38    \\
    800   & 8.00  & 1.73 $\cdot 10^{-1}$ & 0.16  & 5.12 $\cdot 10^{-2}$ & 0.27    \\
    \hline
    \end{tabular}
    \end{center}
    \caption{A comparison of the numerical errors $E$ and CPU times for the implicit first order scheme and the compact implicit second order scheme with the constant reconstruction parameter $\omega=1/2$ for the smooth case of the linear benchmark problem in Experiment~\ref{exp2}, with the exact Gaussian solution given by \eqref{linear_gaussian_solution}. The computations were performed using $M=640$ finite volumes for increasing Courant numbers satisfying $\mathrm{Co}_{L}=100\,\mathrm{Co}_{S}$.}
    \label{TAB:2example_compareerrorwithcn_gaussians_new}
\end{table}

\begin{figure}[H]
    \begin{center}
        \includegraphics[width=0.7\paperwidth]{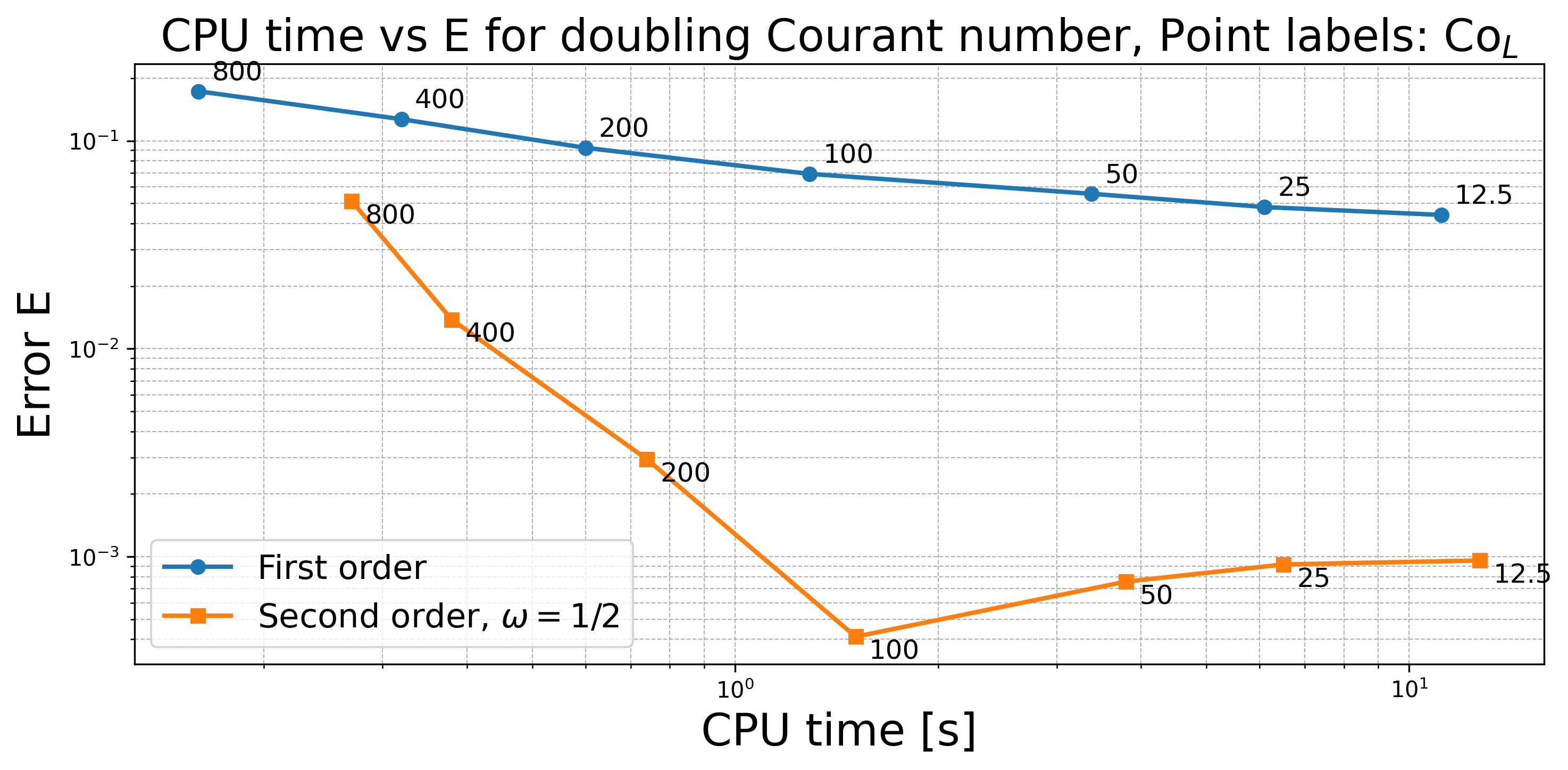}
        \caption{Efficiency comparison of the implicit schemes for the smooth case of the linear benchmark problem in Experiment~\ref{exp2}, with the exact Gaussian solution given by \eqref{linear_gaussian_solution}. The numerical error $E$ is plotted as a function of CPU time for the implicit first order scheme and the compact implicit second order scheme with the constant reconstruction parameter $\omega=1/2$. The plotted points correspond to the larger Courant numbers $\mathrm{Co}_{L}$  reported in Table~\ref{TAB:2example_compareerrorwithcn_gaussians_new}, with $\mathrm{Co}_{L}=100\,\mathrm{Co}_{S}$ and $M=640$ finite volumes.}
        \label{FIG:2example_compareerrorwithcn_gaussians_new}
    \end{center}
\end{figure}

For this experiment, the choice $\mathrm{Co}_{L}=100$, corresponding to $\mathrm{Co}_{S}=1$, gives the best overall performance when the second order scheme with constant $\omega=1/2$ is used. This choice leads to the smallest errors among the tested Courant numbers, while still requiring a reasonable CPU time, see Table~\ref{TAB:2example_compareerrorwithcn_gaussians_new} and also noticeable in Figure~\ref{FIG:2example_compareerrorwithcn_gaussians_new}. This observation is consistent with \cite{frolkovic2023high, zakovafrolkovivc2026esaim}, where different values of $\omega$ were tested and the choice $\omega=1/2$ was shown to give the smallest errors overall. Therefore, $\mathrm{Co}_{L}=100$ provides the most favorable balance between accuracy and computational efficiency for this smooth test case.

For the discontinuous case, we prescribe a rectangular pulse associated with the slow characteristic component, again setting the fast component to zero. Since the slow wave propagates with unit speed, the exact solution is
\begin{equation}
    c_1(x,t)=\begin{cases}
    0.5, & \text{if } 1+t \leq x \leq 2+t\,\\
    0, & \text{otherwise}
    \end{cases}\,, \\ \quad
    c_2(x,t)=\begin{cases}
    -0.5, & \text{if }  1+t \leq x \leq 2+t\,\\
    0, & \text{otherwise}
    \end{cases} \,.
\end{equation}
In contrast to the Gaussian case, this problem requires limiting in order to control oscillations near the discontinuities. We therefore compare several second order implicit variants for increasing Courant numbers (see Table~\ref{TAB:2example_compareerrorwithcn_discontinuous_new}).

\begin{table}[h!]    
\small{
    \begin{center} 
    \begin{tabular}{c c | c c | c c | c c }
    \hline
    \multicolumn{2}{c}{} & \multicolumn{2}{c}{Second order with WENO} & \multicolumn{2}{c}{HR with WENO - limiter with $\mathrm{Co}_{L}$} & \multicolumn{2}{c}{HR with WENO - limiter with $\mathrm{Co}_{S}$} \\
    \hline
    $\mathrm{Co}_{L}$ & $\mathrm{Co}_{S}$ &  $E$   & CPU time [s] & $E$   & CPU time [s] & $E$   & CPU time [s] \\ 
    \hline
    12.5  & 0.125 & 3.99 $\cdot 10^{-2}$ & 36.62 & 5.92 $\cdot 10^{-2}$ & 235.07 & 3.94 $\cdot 10^{-2}$ & 245.54  \\
    25    & 0.25  & 4.32 $\cdot 10^{-2}$ & 19.05 & 7.96 $\cdot 10^{-2}$ & 118.94 & 4.02 $\cdot 10^{-2}$ & 115.44  \\
    50    & 0.50  & 5.08 $\cdot 10^{-2}$ & 8.03  & 1.00 $\cdot 10^{-1}$ & 60.13  & 4.18 $\cdot 10^{-2}$ & 56.96   \\    
    100   & 1.00  & 4.96 $\cdot 10^{-2}$ & 3.54  & 1.26 $\cdot 10^{-1}$ & 30.10  & 4.51 $\cdot 10^{-2}$ & 32.04   \\    
    200   & 2.00  & 6.59 $\cdot 10^{-2}$ & 3.32  & 1.61 $\cdot 10^{-1}$ & 15.03  & 5.16 $\cdot 10^{-2}$ & 15.95   \\
    400   & 4.00  & 1.07 $\cdot 10^{-1}$ & 1.00  & 2.12 $\cdot 10^{-1}$ & 10.66  & 6.95 $\cdot 10^{-2}$ & 8.68    \\
    800   & 8.00  & 2.01 $\cdot 10^{-1}$ & 0.56  & 2.88 $\cdot 10^{-1}$ & 6.37   & 1.16 $\cdot 10^{-1}$ & 3.79    \\
    \hline
    \end{tabular}
    \end{center}
    \caption{A comparison of the numerical errors $E$ and CPU times for the compact implicit second order scheme with WENO reconstruction without time limiting, the compact implicit high-resolution scheme with WENO reconstruction and the time limiter based on $\mathrm{Co}_{L}$, and the same high-resolution scheme with the time limiter based on $\mathrm{Co}_{S}$ for the discontinuous case of the linear benchmark problem in Experiment~\ref{exp2}. The computations were performed using $M=640$ finite volumes for increasing Courant numbers satisfying $\mathrm{Co}_{L}=100\,\mathrm{Co}_{S}$. 
    }
    \label{TAB:2example_compareerrorwithcn_discontinuous_new}}
\end{table}

\begin{figure}[h]
    \begin{center}
        \includegraphics[width=0.7\paperwidth]{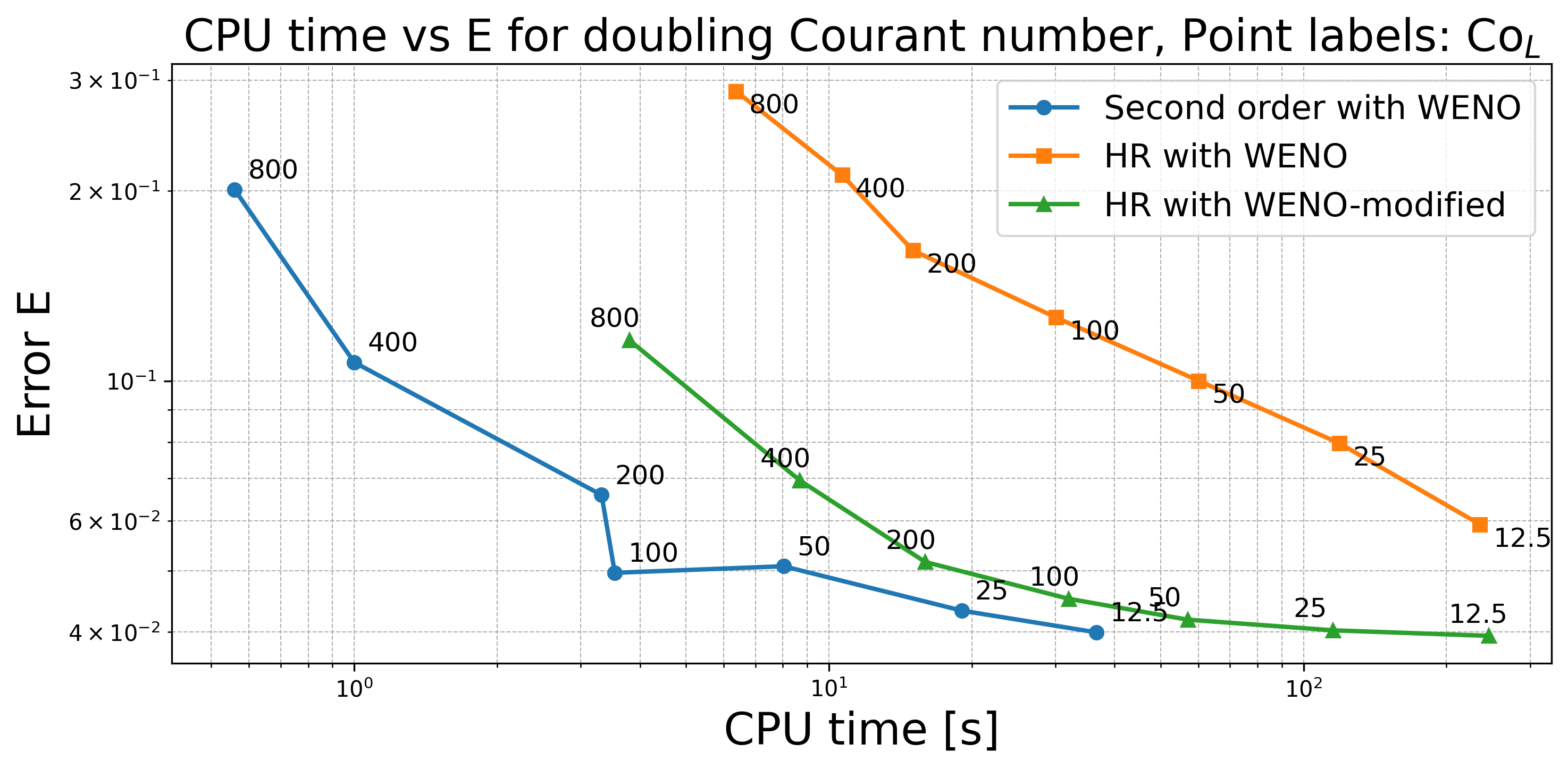}
        \caption{Efficiency comparison of the compact implicit schemes for the discontinuous case of the linear benchmark problem in Experiment~\ref{exp2}, with the rectangular pulse exact solution. The numerical error $E$ is plotted as a function of CPU time for the second order scheme with WENO reconstruction without time limiting, the high-resolution scheme with WENO reconstruction and the time limiter based on $\mathrm{Co}_{L}$, and the high-resolution scheme with WENO reconstruction and the time limiter based on $\mathrm{Co}_{S}$. The plotted points correspond to the Courant numbers $\mathrm{Co}_{L}$ reported in Table~\ref{TAB:2example_compareerrorwithcn_discontinuous_new}, with $\mathrm{Co}_{L}=100\,\mathrm{Co}_{S}$ and $M=640$ finite volumes.}
        \label{FIG:cpu_error_courant_doubling_stvorec_2ndversion}
    \end{center}
\end{figure}

\begin{figure}[H]
    \begin{center}
        \includegraphics[width=0.8\paperwidth]{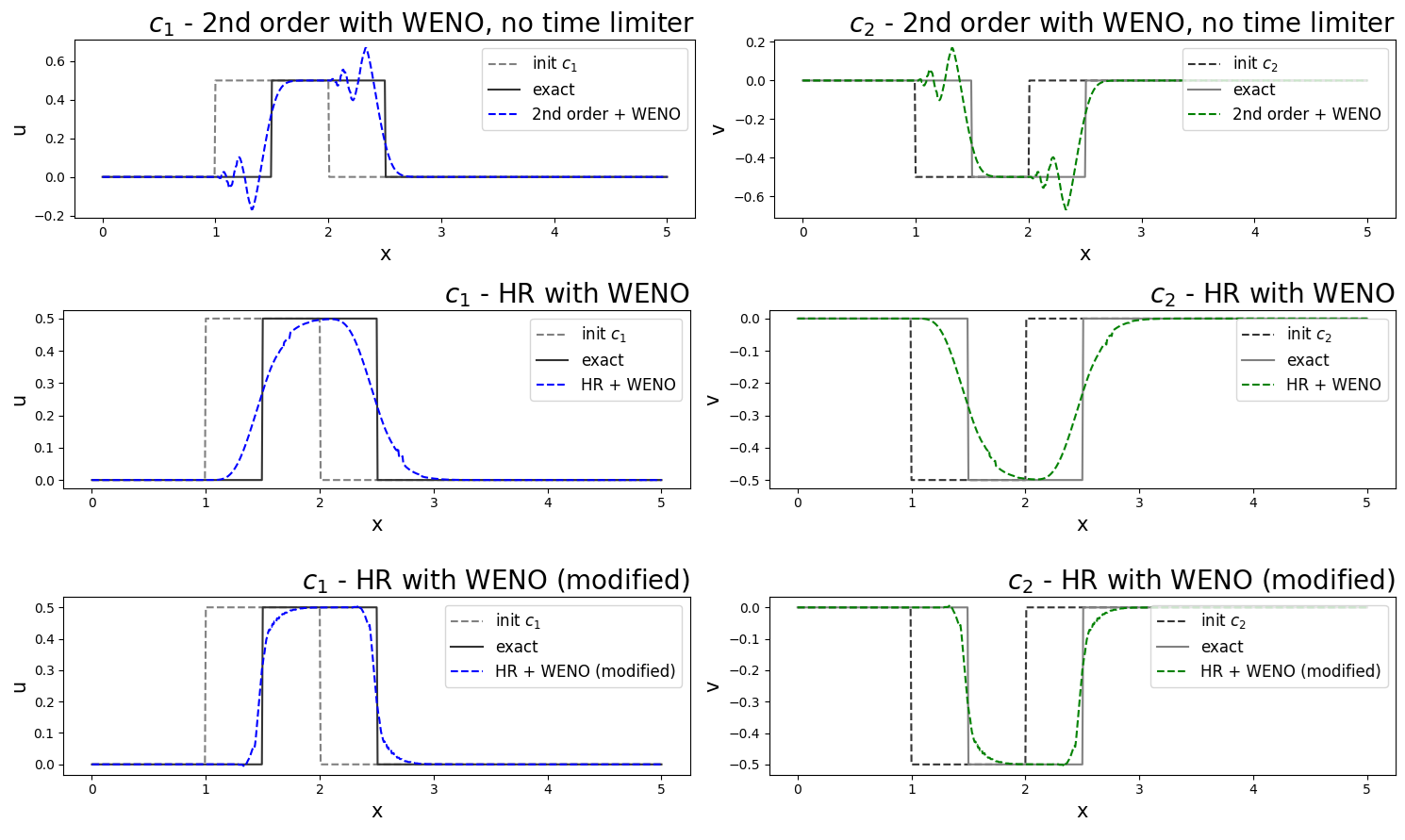}
        \caption{Comparison of the numerical concentration profiles for the discontinuous case of the linear benchmark problem in Experiment~\ref{exp2}, with $\mathrm{Co}_{L}=800$ and $\mathrm{Co}_{S}=8$. The results are shown for the compact implicit second order scheme with WENO reconstruction without time limiting (first row), the compact implicit high-resolution scheme with the time limiter based on $\mathrm{Co}_{L}$ (second row), and the same high-resolution scheme with the time limiter based on $\mathrm{Co}_{S}$ (third row). The WENO scheme without time limiting exhibits visible oscillations near the discontinuities, whereas both high-resolution variants control these oscillations. Using $\mathrm{Co}_{S}$ in the limiter results in sharper and less diffusive concentration profiles than using $\mathrm{Co}_{L}$.}
        \label{FIG:2.priklad_stvorce_porovnanie_cn800}
    \end{center}
\end{figure}

First, we consider the second order with only space limiting (WENO approximation), while the time limiter is omitted. Although the corresponding errors in Table~\ref{TAB:2example_compareerrorwithcn_discontinuous_new} appear promising, the numerical profiles in Figure~\ref{FIG:2.priklad_stvorce_porovnanie_cn800} show visible oscillations for very large Courant numbers, for example $\mathrm{Co}_{L}=800$ and $\mathrm{Co}_{S}=8$. Thus, space limiting alone is not sufficient in this regime.

Second, we apply the full high-resolution scheme WENO approximation, combining both space and time limiting. When the time limiter is based on the maximal Courant number $\mathrm{Co}_{L}$, the method becomes more restrictive and the solution is locally driven closer to first order accuracy. This leads to larger errors, but improves the non-oscillatory behavior of the numerical solution see Figure~\ref{FIG:2.priklad_stvorce_porovnanie_cn800}.

Finally, we test a modified time-limiting strategy, where the limiter is based only on the relevant smaller Courant number $\mathrm{Co}_{S}$ associated with the slow characteristic component present in the solution. Since the fast component is absent in this experiment, this choice reduces unnecessary limiting. As a result, the method produces smaller errors and gives more accurate profiles than the limiter based on $\mathrm{Co}_{L}$, while still controlling the oscillations near the discontinuities (see Figure~\ref{FIG:2.priklad_stvorce_porovnanie_cn800}).

Since the exact solution contains discontinuities, small numerical oscillations may appear near the wave fronts, particularly when higher-order reconstructions are used. Such behavior is a typical feature of high-resolution schemes applied to nonsmooth data \cite{jiang1996efficient, leveque_finite_2004, shu1988efficient, shu_essentially_1998, harten1997high, harten1997uniformly}. In the present computations, however, these oscillations are not unstable oscillations in the sense of a growing numerical instability. They remain bounded, localized near the discontinuities, and do not destroy the solution. The purpose of the weighted essentially non-oscillatory reconstruction is precisely to minimize such spurious oscillations by adapting the reconstruction stencil in nonsmooth regions.

If a strictly non-oscillatory profile is required, the reconstruction may be made more dissipative by applying a stronger limiter, locally reducing the scheme to first order near discontinuities, or using a characteristic-wise reconstruction.

\subsection{Third experiment: three-component displacement chromatography}
\label{subsec:3experiment}

In the last numerical experiment, we consider a three-component chromatography problem \cite{guMathematicalModelingScaleUp2015,
schmidt-traubPreparativeChromatography2020, donat2018implicit} motivated by the test case proposed in \cite{javeed2011efficient, frenz1988high} and later used in \cite{donat2024weno}. The process models the separation of a mixture in a column filled with a solid stationary phase and a liquid moving phase. The mixture contains two solutes, denoted by $c_1$ and $c_2$, while the third component, $c_3$, represents the displacer \cite{frenz1988high,guMathematicalModelingScaleUp2015, schmidt-traubPreparativeChromatography2020, javeed2011efficient}.

The first component is assumed to be less retained than the second one, while the displacer has the strongest affinity to the stationary phase. Therefore, after the displacer is introduced into the column, it progressively replaces the sample components from the stationary phase and pushes them through the column. For sufficiently long columns and suitable adsorption parameters, this process leads to the formation of separated concentration zones, usually referred to as a displacement train or an isotachic train \cite{javeed2011efficient, frenz1988high,donat2024weno, 
mazzotti2009nonclassical, schmidt-traubPreparativeChromatography2020}.

This type of example is commonly used as a benchmark for numerical methods for nonlinear chromatography, since it contains sharp moving transitions and strongly coupled concentration fronts
\cite{javeed2011efficient, donat2018implicit, burger_linearly_2018,donat2024weno, guiochon2006fundamentals}. It tests whether the method can correctly capture the propagation speed of sharp concentration fronts, preserve the separation between components, and reproduce the formation of the displacement train without introducing significant numerical diffusion or spurious oscillations.

Following \cite{javeed2011efficient, donat2024weno, mazzotti2009nonclassical}, we use the Langmuir parameters $a_1=4, \, a_2=5, \, a_3=6,$ and $b_1=4, \, b_2=5, \, b_3=1.$ The porosity parameter is set to $\varepsilon=0.5$ and velocity $u=20$.

In contrast to some computations in \cite{donat2024weno, donat2018implicit, javeed2011efficient, burger_linearly_2018}, where generalized Langmuir-type isotherms are also considered, we restrict ourselves to the classical Langmuir isotherm. Hence, the additional heterogeneity parameter is set to one. Moreover, in the present computations we neglect axial dispersion and set ($D=0$).

Initially, the column contains no transported components, $C(x,0)=(0,0,0)$. The components enter the column through the inlet boundary by means of a prescribed injection profile, and the separation process is generated entirely by the time-dependent inlet data.

The components $c_1$ and $c_2$ are injected only during the initial time interval $0\leq t\leq 0.1$, both with concentration equals $1$. After this injection period, the displacer $c_3$ is injected continuously with concentration equals $1$ for the first simulation ( $a)$ in \eqref{inletconcentrations}), and $c_3=0.5$ for second simulation ( $b)$ in \eqref{inletconcentrations}), the inlet concentrations given by 

\begin{alignat}{2}
\label{inletconcentrations}
\textit{a)}\quad
C_{\mathrm{inj}}(t) &=
\begin{cases}
(1,1,0), & 0\leq t\leq 0.1,\\
(0,0,1), & t>0.1,
\end{cases}
\qquad
&
\textit{b)}\quad
C_{\mathrm{inj}}(t) &=
\begin{cases}
(1,1,0), & 0\leq t\leq 0.1,\\
(0,0,0.5), & t>0.1.
\end{cases}
\end{alignat}

\begin{figure}[h]
    \begin{center}
        \includegraphics[width=0.8\paperwidth]{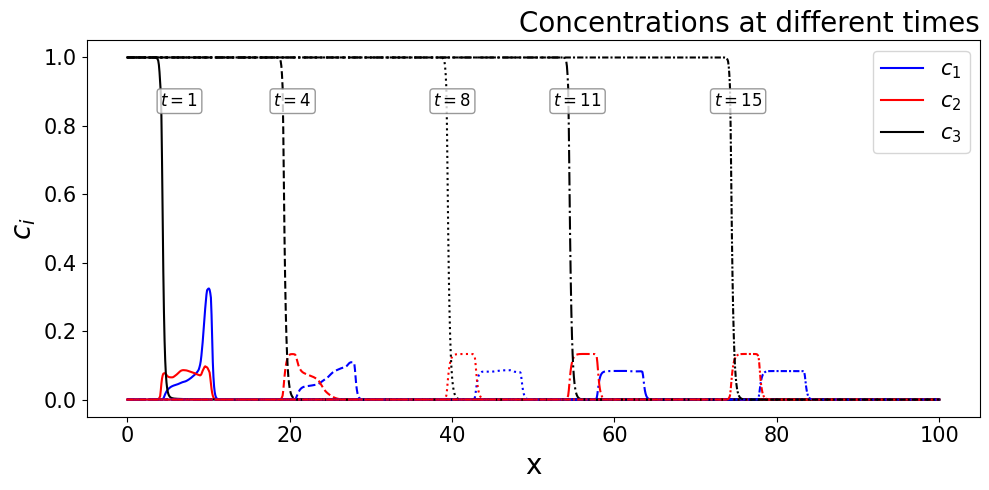}
        \caption{Evolution of the concentration profiles $c_1$ (blue), $c_2$ (red), and $c_3$ (black) for the first inlet profile in \eqref{inletconcentrations}, computed using the high-resolution compact implicit scheme with WENO reconstruction. The solutions are shown at $t\approx 1,\,4,\,8,\,11,$ and $T=15$, with $u=20$ and $\mathrm{Co}_{\max}=\tau u/h=12$.}
        \label{FIG:3experiment_1}
    \end{center}
\end{figure}

In the first simulation, the computational domain is chosen for $x\in[0,100]$ and the mesh contains $M=1280$ finite volumes. The solution is computed up to the final time $T=15$. In Figure~\ref{FIG:3experiment_1}, we show the numerical solution at several times, namely $t=1$, $t=4$, $t=8$, $t=11$, and $T=15$. The computation is performed using the high-resolution scheme with WENO reconstruction.

For this simulation, the maximum Courant number obtained from the global bound \eqref{Courantnumber_max} is $\mathrm{Co}_{\max}=12$. However, this value was not used directly in the computational algorithm, since it would lead to excessive time limiting and would force the method too strongly towards first order accuracy. 
Instead, the relevant local characteristic Courant numbers were computed from the characteristic speeds, see \eqref{Courantnumber_characteristic_components}. In the implementation, this was done using the standard NumPy routines for the inverse matrix and eigenvalue computation.

At the final time step, the maximum values over all finite volumes
for the three characteristic fields were approximately
$3.43$, $4.00$, and $4.80$. The largest of these values is the last mentioned, $4.80$. Although this is significantly smaller than the nominal bound $\mathrm{Co}_{\max}=12$, it is still larger than the stability limit usually required by explicit schemes. This justifies the use of the implicit high-resolution method for this experiment.

The numerical profiles shown in Figure~\ref{FIG:3experiment_1} behave as expected for a displacement chromatography process. After the initial injection of components $c_1$ and $c_2$, the displacer $c_3$ enters the column and progressively pushes the sample components forward. As the process develops, the initially overlapping concentrations separate into distinct zones with increasingly sharp interfaces. At the final time, the three components form a well-developed sequence of concentration profiles, indicating an effective displacement and separation process. This behavior is consistent with the corresponding displacement chromatography experiment reported in \cite{javeed2011efficient,donat2024weno}.

\begin{figure}[H]
    \begin{center}
        \includegraphics[width=0.8\paperwidth]{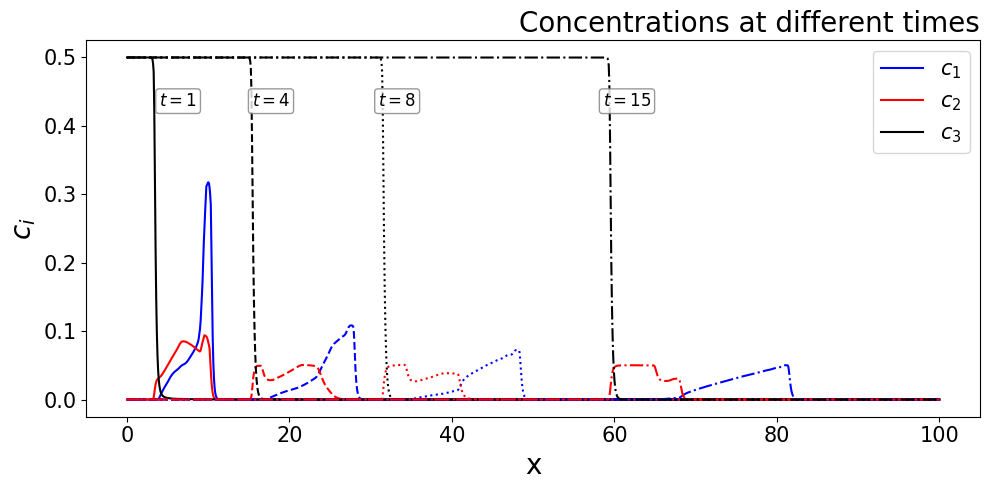}
        \caption{Evolution of the concentration profiles $c_1$ (blue), $c_2$ (red), and $c_3$ (black) for the second inlet profile in \eqref{inletconcentrations}, with the displacer concentration set to $c_3=0.5$. The solutions are computed using the high-resolution compact implicit scheme with WENO reconstruction and are shown at $t\approx 1,\,4,\,8,$ and $T=15$, with $u=20$ and $\mathrm{Co}_{\max}=\tau u/h=12$.}
        \label{FIG:3experiment_2}
    \end{center}
\end{figure}

In the second simulation, we use the same computational domain, spatial mesh, final time, and numerical method as in the first simulation. The only difference is the lower concentration of the continuously injected displacer, which is reduced from $c_3=1$ to $c_3=0.5$ according to the second inlet profile in \eqref{inletconcentrations}. The evolution of the three concentration profiles is shown in Figure~\ref{FIG:3experiment_2}.

Compared with the first simulation, the lower displacer concentration produces a weaker displacement of the initially injected components. Consequently, the separation of $c_1$ and $c_2$ develops more slowly and remains less pronounced at the final time. In particular, the first component does not form the same distinct rectangular concentration zone as in the case $c_3=1$, while the profile of the second component exhibits a less developed staircase-like structure. Thus, reducing the amount of injected displacer prevents the complete formation of the isotachic train within the considered time interval. This qualitative behavior is consistent with the observations reported in \cite{donat2024weno} for a reduced displacer concentration.

Overall, the compact implicit high-resolution scheme reproduced the expected displacement behavior with well-resolved concentration fronts and without any visible oscillations, even for the relatively large Courant number $\mathrm{Co}_{\max}=12$.

\section{Conclusions}

In this paper, we applied compact implicit high-resolution schemes to multicomponent transport problems with Langmuir sorption. The main difficulty of the system case is the nonlinear coupling introduced by the multicomponent Langmuir isotherm, which leads to a nonlinear retardation matrix and hyperbolic system with several characteristic speeds.

The numerical experiments show that the proposed compact implicit schemes achieve the expected accuracy for smooth solutions and remain stable for Courant numbers significantly larger than those admissible for explicit schemes, which still require the nonlinear inversion of the accumulation function and are restricted by a much stronger stability condition.

Special attention was paid to the Courant number used in the limiter. The use of the global bound $\mathrm{Co}_{\max}=\tau u/h$ may lead to excessive limiting, especially when the actual characteristic speeds are reduced by sorption. Therefore, characteristic Courant numbers based on the eigenvalues of $J(C)^{-1}$ were used in selected computations. This strategy reduces unnecessary numerical diffusion and allows the high-resolution scheme to preserve more of its second order accuracy.

The chromatography examples demonstrate that the method is able to capture sharp moving fronts and the formation of displacement trains in multicomponent Langmuir systems.

\bibliographystyle{plain}
\bibliography{lit.bib}

\end{document}